\documentclass[11pt]{article}
\usepackage{graphicx}
\usepackage{amsmath, amsthm, amssymb}
\usepackage{epstopdf}
\usepackage[plainpages=false,pdfpagelabels,bookmarksnumbered]{hyperref}
\newtheoremstyle{thm}{\topsep}{\topsep}
     {}
     {}
     {\bfseries}
     {}
     {\newline}
     {\thmname{#1}\thmnumber{ #2}\thmnote{ #3}}

   \theoremstyle{thm}
   \newtheorem{thm}{Theorem}[section]

\newtheoremstyle{cor}{\topsep}{\topsep}
     {}
     {}
     {\bfseries}
     {}
     {\newline}
     {\thmname{#1}\thmnumber{ #2}\thmnote{ #3}}

   \theoremstyle{cor}
   \newtheorem{cor}[thm]{Corollory}

\newtheoremstyle{lem}{\topsep}{\topsep}
     {}
     {}
     {\bfseries}
     {}
     {\newline}
     {\thmname{#1}\thmnumber{ #2}\thmnote{ #3}}

   \theoremstyle{lem}
   \newtheorem{lem}[thm]{Lemma}

\newtheoremstyle{prop}{\topsep}{\topsep}%
     {}
     {}
     {\bfseries}
     {}
     {\newline}
     {\thmname{#1}\thmnumber{ #2}\thmnote{ #3}}

   \theoremstyle{prop}

\newtheoremstyle{remark}{\topsep}{\topsep}%
     {}
     {}
     {\bfseries}
     {}
     {\newline}
     {\thmname{#1}\thmnumber{ #2}\thmnote{ #3}}

   \theoremstyle{remark}
   \newtheorem{rem}[thm]{Remark}
   
\newtheoremstyle{defintion}{\topsep}{\topsep}%
     {}
     {}
     {\bfseries}
     {}
     {\newline}
     {\thmname{#1}\thmnumber{ #2}\thmnote{ #3}}

   \theoremstyle{defintion}
   \newtheorem{de}[thm]{Defintion}

\newtheorem*{thm4.1}{Theorem 4.5}
\newtheorem*{thm5.1}{Theorem 5.3}

\numberwithin{equation}{section}

\def \Z{\mathbb Z}
\def \C{\mathbb C}

\def \N{\mathbb N}

\def \wt{{\rm wt}}

\def \span{{\rm span}}

\def \wt{{\rm wt}}
\def \deg{{\rm deg}}

\def \Res{{\rm Res}}

\def \End{{\rm End}}
\def \Hom{{\rm Hom}}

\def \<{\langle}
\def \>{\rangle}

\def \pf{\noindent {\bf Proof:} \,}
\def \cg{\chi_g}
\def \cg'{\chi'_g}

\def \d{\delta}

\def \r{\rho}

\def \1{{\textbf 1}}

\def \qed{{\mbox{$\square$}}}

\begin{document}

\begin{center}
{\Large {\bf Ordered spanning sets for quasimodules for M\"{o}bius vertex algebras}}
\\
\vspace{0.5cm} Geoffrey Buhl\footnote{Supported by a California State University Channel Islands Faculty Development Grant.}\\
Mathematics Program\\
California State University Channel Islands\\
 Camarillo, CA 93012\\
{\tt geoffrey.buhl@csuci.edu}
\end{center}

\begin{abstract}
Quasimodules for vertex algebras are generalizations of modules for vertex algebras.  These new objects arise from a generalization of locality for fields.  Quasimodules tie together module theory and twisted module theory, and both twisted and untwisted modules feature Poincar\'{e}-Birkhoff-Witt-like spanning sets.  This paper generalizes these spanning set results to quasimodules for certain M\"{o}bius vertex algebras.  In particular this paper presents two spanning sets, one featuring a difference-zero ordering restriction on modes and another featuring a difference-one ordering restriction.
\end{abstract}

\section{Introduction}

Representation theory is a particularly rich aspect of the theory of vertex algebras with fundamental connections to number theory, the theory of simple groups, and string and conformal field theories in physics.  Quasimodules for vertex algebras are new module-like structures whose vertex algebra action is governed by a modified Jacobi identity.  These new objects are generalizations of modules for vertex algebras,  and quasimodules are related to twisted modules for vertex algebras.  There are Poincar\'{e}-Birkhoff-Witt-like spanning sets for both modules and twisted modules for vertex operator algebras.  This paper extends these spanning set results to quasimodules for the more general M\"{o}bius vertex algebras.

The theory of quasimodules for vertex algebras was developed by Li \cite{MR1966260} \cite{MR2215259} \cite{MR2218823} \cite{MR2298863} \cite{li-2006}.  Quasimodules for vertex algebras arise from a natural generalization of locality.  For a vector space $W$, two elements $a(x), b(x)$ of $\Hom(W,W((x)))$, called weak vertex operators, are {\it local} if there exists a non-negative integer $k$ such that 
\begin{eqnarray}
(x_1-x_2)^k a(x_1)b(x_2)=(x_1-x_2)^k b(x_2)a(x_1).
\end{eqnarray}
Maximal local subsets of $\Hom(W,W((x)))$ are vertex algebras with $W$ as a module \cite{MR1387738}.  Quasimodules arise from a weaker form of locality called quasi-locality.  Two weak vertex operators $a(x), b(x) \in \Hom(W,W((x)))$ are {\it quasi-local} if there exists a non-zero polynomial $f(x_1,x_2) \in \C[x_1,x_2]$ such that 
\begin{eqnarray}
f(x_1,x_2) a(x_1)b(x_2)=f(x_1,x_2)b(x_2)a(x_1).
\end{eqnarray}
One surprising result is: maximal quasi-local subsets of $\Hom(W,W((x)))$ are vertex algebras  \cite{MR2218823}.  That is, quasi-locality is equivalent to locality in the construction of vertex algebras, and as the Jacobi identity axiom can be replaced by a locality axiom in the definition of a vertex algebra, the Jacobi identity axiom can be replaced by a quasi-locality axiom in the definition of a vertex algebra. This construction of vertex algebras from quasi-local fields does not also result in modules for the vertex algebras as in the locality construction.  This quasi-locality construction creates module-like structures called quasimodules, and the main axiom for these new objects is a modified version of Jacobi identity for vertex algebras and their modules.

Quasimodule is a new `super-type' for modules for vertex algebras and as a category contain as a subcatogory modules for vertex algebras.  Li shows that there is a natural isomorphism between twisted modules for vertex operator algebras and quasimodules for certain vertex operator algebras \cite{li-2006}.  As such the theory of quasimodules encompasses the theory modules and twisted-modules for vertex (operator) algebras.  In this paper, we extend a result common to the theory of module and twisted-modules to the the theory of quasimodules.  

There are a number of Poincar\'{e}-Birkhoff-Witt-like spanning sets for vertex (operator) algebras, modules, and twisted modules. \cite{MR1700507} \cite{MR1990879} \cite{MR1927435} \cite{MR2046807} \cite{MR2039213}.  In each spanning set, elements have the form
\begin{eqnarray}
u^{(1)}_{n_1}\cdots u^{(k)}_{n_k}w,
\end{eqnarray}
where $w$ is the vacuum vector for algebra spanning sets or a generating vector for module spanning sets. Here the $u^{(i)}$'s are elements of some subspace of the vertex (operator) algebra, and there are ordering restrictions on the indexes $n_i$. These restrictions appear as difference conditions, similar to difference conditions on partitions.  A difference-$n$ condition on indicies means that the indices of adjacent modes must differ by at least $n$. That is, for adjacent modes $u^{(i)}_{m_i}$ and $u^{(i+1)}_{m_{i+1}}$ in a spanning set element, $|m_i-m_{i+1}| \geq n$.

Work by Gaberdiel and Neitzke, Miyamoto, Yamauchi, and the author has shown that there exists difference-one spanning sets for both modules and twisted modules for vertex operator algebras, and work of Karel and Li shows that there exist difference-zero algebra and module spanning sets.  The main results presented in this paper are analogous results for quasimodules.  In particular there exists difference-zero and  difference-one spanning sets for quasimodules for $\N$-graded M\"{o}bius vertex algebras.

\begin{thm4.1}
For a M\"{o}bius vertex algebra $V$, $X$ a set of homogeneous representatives of a basis for $V/C_1(V)$, and a quasimodule module $W$ generated by a vector $w$ that is uniformly annihilated by $X$, $W$ is spanned by elements of the form
\begin{eqnarray}
x^{(1)}_{n_1} \cdots x^{(r)}_{n_r}w
\end{eqnarray}
with $n_1\leq \cdots \leq n_r<T$ where $r \in \N$, $x^{(1)}, \ldots, x^{(r)} \in X$, $n_1, \ldots , n_r \in \Z$, and $T$ is order of uniform annihilation of $w$ by $X$.
\end{thm4.1}

\begin{thm5.1}
For a M\"{o}bius vertex algebra $V$, $X$ a set of homogeneous representatives of a basis for $V/C_2(V)$, and a quasimodule module $W$ generated by a vector $w$ that is uniformly annihilated by $X$, $W$ is spanned by elements of the form
\begin{eqnarray}
x^{(1)}_{n_1} \cdots x^{(r)}_{n_r}w
\end{eqnarray}
with $n_1< \cdots < n_r<T$ where $r \in \N$, $x^{(1)}, \ldots, x^{(r)} \in X$, $n_1, \ldots , n_r \in \Z$, and $T$ is order of uniform annihilation of $w$ by $X$.
\end{thm5.1}

Depending on one's point of view, M\"{o}bius vertex algebras are either generalizations of vertex operator algebras and conformal vertex algebras or restrictions of vertex algebras.  Vertex operator algebras feature a representation of the Virasoro algebra and a lower truncation assumption on their grading.  Conformal vertex algebras admit a Virasoro representation, but do not presuppose a lower truncation assumption for their grading.    M\"{o}bius vertex algebras admit a representation of $\mathfrak{sl}(2)$ rather than the Virasoro algebra, and no lower-truncation on their grading is assumed.  Vertex algebras are more general objects with no grading and no Lie algebra representation assumed.  The  main results of this paper rely on the grading of M\"{o}bius vertex algebras, and in particular some results require a lower-truncated grading.  Some prior spanning sets required vertex operator algebras that are $C_2$-cofinite or of CFT-type.  The results of this paper apply in a more general setting.

The outline of this paper is as follows: the second section reviews M\"{o}bius vertex algebras and modules for M\"{o}bius vertex algebras.  Section three presents the theory of quasimodules and derives a number of identities for quasimodules from the quasi-Jacobi identity.  The third section also presents `quasi' analogues to the associativity and commutativity identities.  The fourth section contains the proof of a difference-zero spanning set for quasimodules.  Section five contains the proof of the  difference-one spanning set for quasimodules and results relating $C_2$-cofiniteness and $C_n$-cofiniteness for quasimodules for certain M\"{o}bius vertex algebras.

\section{ M\"{o}bius vertex algebras and modules}

In this section we present the definition of a M\"{o}bius vertex algebra and modules for these algebras as well as some related concepts.  For an introduction to the theory of vertex algebras and vertex operator algebras, we refer the reader to \cite{MR2023933}.

M\"{o}bius vertex algebras appear in the work of Huang, Lepowsky, and Zhang on a logarithmic tensor product theory for conformal vertex algebras \cite{huang-2006}.  These algebras are generalized versions of conformal vertex algebras, $\Z$-graded vertex algebras that admit a Virasoro algebra representation.  M\"{o}bius vertex algebras only admit a representation of $\mathfrak{sl}(2)$, a Lie subalgebra of the Virasoro algebra.

\begin{de}
A {\it M\"{o}bius vertex algebra} is a $\Z$-graded vector space
 \begin{eqnarray}
V=\coprod_{n \in \Z} V_n
\end{eqnarray}
equipped with a linear map 
\begin{eqnarray}
Y:  V &\rightarrow& \End(V)[[x,x^{-1}]]\\
  v &\mapsto& Y(v,x)=\sum_{n \in \Z}v_n x^{-n-1} \ \mbox{(where $V_n \in \End(V)$)}
\end{eqnarray}
where $Y(v,x)$ is called the vertex operator associated with $v$ and a distinguished vector $\1 \in V_0$ (the vacuum vector), satisfying the following conditions for $u,v \in V$: the lower truncation condition: 
\begin{eqnarray}
u_nv=0 \mbox{ for n sufficiently large};
\end{eqnarray}
the vacuum property:
\begin{eqnarray}
Y(\1,x)=1_V;
\end{eqnarray}
the creation property:
\begin{eqnarray}
Y(v,x) \1 \in V[[x]] \mbox{ and } \lim_{x \rightarrow 0} Y(v,x) \1=v;
\end{eqnarray}
the Jacobi Identity:
\begin{eqnarray}
&x_0^{-1}\d \left(\frac{x_1 - x_2}{x_0}\right)Y(u,x_1)Y(v,x_2)-
x_0^{-1} \d \left(\frac{x_2- x_1}{-x_0}\right)Y(v,x_2)Y(u,x_1) \nonumber \\
&=x_2^{-1} \d \left(\frac{x_1- x_0}{
x_2}\right)Y(Y(u,x_0)v,x_2).
\end{eqnarray}
In addition there is a representation $\rho$ of $\mathfrak{sl}(2)$ on V given by:
\begin{eqnarray}
L(j)=\rho(L_j) \mbox{, } j=-1,0,1
\end{eqnarray}
where $\{ L_{-1}, L_0, L_1\}$ from a basis of $\mathfrak{sl}(2)$ with Lie brackets
\begin{eqnarray}
[L_0,L_{-1}]= L_{-1}, [L_0,L_1]=-L_1  \mbox{, and} [L_{-1},L_{1}]=-2L_0, 
\end{eqnarray}
and the following conditions hold for $v \in V$ and $j=-1,0,1$:
\begin{eqnarray}
[L(j), Y(u,x)]= \sum^{j+1}_{k=0} \binom{j+1}{k}x^{j+1-k}Y(L(k-1)v,x),
\end{eqnarray}

\begin{eqnarray}
\frac{d}{dx}Y(v,x)=Y(L(-1)v,x),
\end{eqnarray}
and
\begin{eqnarray}
L(0)v=nv=(\wt v)v \mbox{ for } n \in \Z \mbox{ and } v \in V_n.
\end{eqnarray}
\end{de}

A M\"{o}bius vertex algebra is denoted by the quadruple $(V,Y,\1,\r)$ or when clear from the context as $V$.

\begin{de}
A M\"{o}bius vertex algebra $V=\coprod_{n \in \Z} V_n$ is {\it $\N$-graded} if $V_n=0$ for $n<0$.
\end{de}

In order to prove spanning set theorems, we must assume that there are no vectors of negative weight in the M\"{o}bius vertex algebra.  However, the identities used to construct the spanning sets hold for any M\"{o}bius vertex algebra.    Specifically the results in the third section do not require the M\"{o}bius vertex algebras to be $\N$-graded, while the major theorems of the fourth and fifth sections do require this assumption.

The definition of a M\"{o}bius vertex algebra is similar to the definition of a quasi-vertex operator algebra given by Frenkel, Huang, and Lepowsky in \cite{MR1142494}.  The definition of a quasi-vertex operator algebra includes two axioms in addition to those of a M\"{o}bius vertex algebra: each graded piece is finite dimensional and the $\Z$-grading is truncated from below.  Though the  $\N$-graded M\"{o}bius vertex algebras that we consider in this paper have the lower truncation property of quasi-vertex operator algebras, we do not assume that the graded pieces of the algebra are finite-dimensional.  

There is an additional benefit of using the more general $\N$-graded M\"{o}bius vertex algebras rather than quasi-vertex operator algebras. It prevents titling this paper ``Ordered spanning sets for quasimodules for quasi-vertex operator algebras''.  Note that quasi-vertex operator algebras and quasimodules for vertex algebra are fundementally different notions arising from two distinct generalizations.  Quasimodules arise from a generalization of locality, and quasi-vertex operator algebras arise from a generalization of the Lie algebra representation the algebra admits ($\mathfrak{sl}(2)$ rather than the Virasoro algebra).


In this paper we construct Poincar\'{e}-Birkhoff-Witt-like spanning sets for quasimodules for $\N$-graded M\"{o}bius vertex algebras.  These spanning sets are analogous to spanning sets for modules and twisted modules for vertex operator algebras \cite{MR1927435} \cite{MR2046807} \cite{MR2039213}.  In turn these module spanning sets are generalizations of vertex operator algebra spanning sets \cite{MR1700507}  \cite{MR1990879}.

\begin{de}
A {\it module} for a  M\"{o}bius vertex algebra $V$ is a vector space $W$ equipped with a linear map
\begin{eqnarray}
Y_W:  &&V \rightarrow \End(W)[[x,x^{-1}]]\\
&&  v \mapsto Y_W(v,x)=\sum_{n \in \Z}v_n x^{-n-1} \text{ where } v_n \in \End V
\end{eqnarray}
satisfying the following conditions:\\
the lower truncation condition: for $ v \in V$ and $w \in W$,
\begin{eqnarray}
v_nw=0 \text{ for $n$ sufficiently large};
\end{eqnarray}
the vacuum property:
\begin{eqnarray}
Y_W(\1,x)=1_W;
\end{eqnarray}
the Jacobi identity: for $u,v \in V$, 
\begin{eqnarray}
&x_0^{-1}\d \left(\frac{x_1 - x_2}{x_0}\right)Y_W(u,x_1)Y_W(v,x_2)-
x_0^{-1} \d \left(\frac{x_2- x_1}{-x_0}\right)Y_W(v,x_2)Y_W(u,x_1) \nonumber \\
& =x_2^{-1} \d \left(\frac{x_1- x_0}{
x_2}\right)Y_W(Y(u,x_0)v,x_2).
\end{eqnarray}
In addition there is a representation $\rho$ of $\mathfrak{sl}(2)$ on W given by:
\begin{eqnarray}
L_W(j)=\rho(L_j) \mbox{, } j=-1,0,1
\end{eqnarray}
where $\{ L_{-1}, L_0, L_1\}$ from a basis of $\mathfrak{sl}(2)$ with Lie brackets
\begin{eqnarray}
[L_0,L_{-1}]= L_{-1}, [L_0,L_1]=-L_1  \mbox{, and} [L_{-1},L_{1}]=-2L_0, 
\end{eqnarray}
and the following conditions hold for $v \in V$ and $j=-1,0,1$:
\begin{eqnarray}
[L_W(j), Y_W(u,x)]= \sum^{j+1}_{k=0} \binom{j+1}{k}x^{j+1-k}Y_W(L(k-1)v,x) \label{modsl2}
\end{eqnarray}
and
\begin{eqnarray}
\frac{d}{dx}Y(v,x)=Y_W(L(-1)v,x). \label{modder}
\end{eqnarray}
\end{de}

A module for a M\"{o}bius vertex algebra is denoted by the pair $(W, Y_W)$ or as $W$.  

\begin{rem}
In the definition of a module for a M\"{o}bius vertex algebra, condition  (\ref{modder}) is actually a consequence of the Jacobi identity and the $\mathfrak{sl}(2)$ representation of the M\"{o}bius vertex algebra. It is included in the definition for clarity.
\end{rem}



M\"{o}bius vertex algebras and their modules are infinite-dimensional except for the trivial case, and a size restriction on these objects is often a useful property.  Of particular importance is algebra $C_2$-cofiniteness, defined below, where a M\"{o}bius vertex algebra $V$ is considered as a module for itself.  Many vertex operator algebras are $C_2$-cofinite; vertex operator algebras constructed from lattices, Kac-Moody Lie algebras, and the Virasoro algebra all share this property.  Assuming $C_2$-cofiniteness of a vertex operator algebra is necessary to prove the modularity of certain trace functions for certain vertex operator algebras \cite{MR1317233} and vertex operator superalgebras \cite{MR2175996}.  $C_2$-cofiniteness of the algebra is sometimes referred to as Zhu's finiteness condition.  Here $C_2$-cofiniteness is generalized to $C_n$-cofiniteness for $n \geq 2$.

\begin{de}
For  $n \geq 2$, a module $W$ for a M\"{o}bius vertex algebra $V$ is {\it $C_n$-cofinite}  if $W/C_n(W)$ is finite-dimensional where
\begin{eqnarray}
C_n(W)=\span\{v_{-n}w: v \in V \mbox{ and } w \in W\}.
\end{eqnarray}
\end{de}

\begin{rem}
The $L(-1)$-derivative property ensures that $C_n(W) \subseteq C_m(W)$ for $n \geq m\geq 2$. As a result $C_n$-cofiniteness of $W$ implies $C_m$-cofiniteness of $W$ for $n \geq m$.  The converse, $C_m$-cofiniteness implies $C_n$-cofiniteness for $n \geq m\geq 2$, holds for certain  vertex algebras \cite{MR1990879} and modules \cite{MR1927435}.
\end{rem}

\begin{rem}
There are different ways to extend the definition of $C_n$-cofiniteness to include $n=1$. This paper uses the following conventions: 
a M\"{o}bius vertex algebra $V$ is {\it $C_1$-cofinite} if $V/C_1(V)$ is finite-dimensional where
\begin{eqnarray}
C_1(V)=\span\{u_{-1}v, L(-1)w: u,v \in \coprod_{n >0} V_n, w \in V \},
\end{eqnarray}
and a module $W$ for a M\"{o}bius vertex algebra $V$  is {\it $c_1$-cofinite} if $W/c_1(W)$ is finite-dimensional where
\begin{eqnarray}
c_1(W)=\span\{u_{-1}v: v \in V \mbox{ and } w \in W \}.
\end{eqnarray}

Algebra $C_1$-cofiniteness in the above form appears in the work of Li \cite{MR1676852} and Karel and Li \cite{MR1700507} on minimal generating sets for spanning sets for vertex operator algebras. Module $c_1$-cofiniteness appears in the work of Nahm \cite{MR1305167}, who studied vertex operator algebras for which certain irreducible modules are $c_1$-cofinite.  Adding the possible confusion around the prefix `quasi', Nahm called such algebras {\it quasirational}.  Quasirationality or $c_1$-cofiniteness of certain irreducible modules is a important assumption in Huang's work on modular tensor categories and the Verlinde conjecture \cite{MR2029793} \cite{MR2151865} \cite{MR2140309}.  Some of Huang's work also requires that the algebras be $C_2$-cofinite, which implies $c_1$-cofiniteness of the modules \cite{MR2052955}. 
\end{rem}

\begin{rem}
If we assume that $V_0=\C \1$ for an $\N$-graded M\"{o}bius vertex algebra, then the $L(-1)$-derivative property ensures that $C_n(V) \subseteq C_1(V)$ for $n \geq2$.  This means that $C_n$-cofiniteness of the algebra implies $C_1$-cofiniteness of the algebra.  A similar statement holds for $c_1$-cofiniteness.  However the converse is not true; the Heisenberg vertex operator algebra is $C_1$-cofinite, but is not $C_2$-cofinite. For algebras $c_1$-cofiniteness is uninteresting as the creation axiom ensures that all vertex algebras are $c_1$-cofinite.
\end{rem}

\begin{rem}
A vertex algebra is said to be `of CFT-type' if $V_0=\C \1$ and $V$ is $\N$-graded.  If $V$ is not of CFT-type then $C_2(V)$ is not necessarily contained in $C_1(V)$ as $V_0 \ne \C \1$. In this case $u_{-2}v$ with $v \in V_0$ and $v \not \in \C\1$ is an element of $C_{2}(V)$ that is not necessarily in $C_1(V)$.  
\end{rem}

Representatives of the  quotient spaces $V/C_1(V)$ and $V/C_2(V)$ are used to generate vertex operator algebra spanning sets \cite{MR1700507} \cite{MR1990879} and module spanning sets \cite{MR1927435} \cite{MR2046807} \cite{MR2039213}.  These results do not rely on the full Virasoro algebra representation, and they hold under the weaker assumption of a representation of $\mathfrak{sl}(2)$ on the algebra.  In particular, these Poincar\'{e}-Birkhoff-Witt-like spanning sets exist for M\"{o}bius vertex algebras and modules for M\"{o}bius vertex algebras.  In this paper we extend these results to quasimodules for $\N$-graded M\"{o}bius vertex algebras.

\section{Quasimodules for vertex algebras}

In this section we present the concept of a quasimodule which results naturally from a generalization of locality.  Also we derive `quasi' analogues to the associativity and commutativity identities for modes of vertex operators for vertex algebras and modules.  From these new identities we derive further identities that are used to prove the difference-zero and difference-one spanning sets for quasimodules of $\N$-graded M\"{o}bius vertex algebras.

A standard construction of local vertex algebras is due to Li \cite{MR1387738}.  For a vector space $W$, any local subalgebra of $\Hom (W, W((x)))$ is a vertex algebra.  Maximally local sets are subalgebras and hence are vertex algebras, and they act faithfully on $W$.  In this way $W$ is a module the vertex algebras constructed in this manner.  Locality is a feature a number of vertex algebras, but there exist examples of non-local vertex algebras \cite{MR1966260} \cite{MR2215259}.


\begin{de}
For a vectorspace $W$, a set $U \in \Hom(W,W((x)))$ is {\it local} if for any $a(x), b(x) \in U$, there exists a non-negative integer $k$ such that
\begin{eqnarray}
(x_1-x_2)^k a(x_1)b(x_2)= (x_1-x_2)^kb(x_2)a(x_1).
\end{eqnarray}
\end{de}

In the definition of a vertex algebra, the Jacobi identity is equivalent to the assumption that vertex operators associated with elements of a vertex algebra are local.  Such a reformulation, replacing the Jacobi identity axiom with a locality axiom, does not work for modules.  Modules lack the creation axiom, and as a result the Jacobi identity for modules is not equivalent to locality.  

In \cite{MR2218823} Li explores the structural implications of a generalization of locality called quasi-locality. Quasi-locality is obtained by replacing $(x_1-x_2)^k$ in the definition of locality by a non-zero polynomial $f(x_1,x_2)$.

\begin{de}
For a vectorspace $W$, a set $U \in \Hom(W,W((x)))$ is {\it quasi-local} if for any $a(x), b(x) \in U$, there exists a non-zero $f(x_1,x_2)\in \mathbb{C}[x_1,x_2]$ such that
\begin{eqnarray}
f(x_1,x_2)a(x_1)b(x_2)= f(x_1,x_2)b(x_2)a(x_1).
\end{eqnarray}
\end{de}

Surprisingly, maximal quasilocal subesets of $\Hom(W,W((x)))$ have the structure of a vertex algebra \cite{MR2218823}. That is, quasi-locality is equivalent to the Jacobi identity in the definition of a vertex algebra, just as locality is.  This new construction of vertex algebras from maximal quasi-local subsets of $\Hom(W,W((x)))$ differs for the previous construction in a significant way: $W$ is not a module for these vertex algebras in the traditional sense.  $W$ has the structure of a quasimodule and is governed by a new Jacobi-like identity called the quasi-Jacobi identity. These new objects are also related to twisted modules for vertex algebras \cite{li-2006}.

\begin{de}
A {\it quasimodule} for a M\"{o}bius vertex algebra $V$ is a vector space $W$ equipped with a linear map
\begin{eqnarray}
Y_W:  &&V \rightarrow \End(W)[[x,x^{-1}]]\\
&&  v \mapsto Y_W(v,x)=\sum_{n \in \Z}v_n x^{-n-1} \text{ where } v_n \in \End W
\end{eqnarray}
satisfying the following conditions:\\
the lower truncation condition: for $ v \in V$ and $w \in W$,
\begin{eqnarray}
v_nw=0 \text{ for $n$ sufficiently large};
\end{eqnarray}
the vacuum property:
\begin{eqnarray}
Y_W(\1,x)=1_W;
\end{eqnarray}
the quasi-Jacobi identity: for $u,v \in V$ there exists a non-zero polynomial $f(x_1,x_2) \in \C[x_1,x_2]$ such that, 
\begin{eqnarray}
&x_0^{-1}\d \left(\frac{x_1 - x_2}{x_0}\right)f(x_1,x_2)Y_W(u,x_1)Y_W(v,x_2)-
x_0^{-1} \d \left(\frac{x_2- x_1}{-x_0}\right)f(x_1,x_2)Y_W(v,x_2)Y_W(u,x_1) 
\nonumber \\
&=x_2^{-1} \d \left(\frac{x_1- x_0}{x_2}\right)f(x_1,x_2)Y_W(Y(u,x_0)v,x_2). \label{qJi}
\end{eqnarray}
In addition there is a representation $\rho$ of $\mathfrak{sl}(2)$ on W given by:
\begin{eqnarray}
L_W(j)=\rho(L_j) \mbox{, } j=-1,0,1
\end{eqnarray}
where $\{ L_{-1}, L_0, L_1\}$ from a basis of $\mathfrak{sl}(2)$ with Lie brackets
\begin{eqnarray}
[L_0,L_{-1}]= L_{-1}, [L_0,L_1]=-L_1  \mbox{, and} [L_{-1},L_{1}]=-2L_0, 
\end{eqnarray}
and the following conditions hold for $v \in V$ and $j=-1,0,1$:
\begin{eqnarray}
[L_W(j), Y_W(u,x)]= \sum^{j+1}_{k=0} \binom{j+1}{k}x^{j+1-k}Y_W(L(k-1)v,x) 
\end{eqnarray}
and
\begin{eqnarray}
\frac{d}{dx}Y(v,x)=Y_W(L(-1)v,x).
\end{eqnarray}
\end{de}

A quasimodule is denoted by the pair $(W, Y_W)$ or  as $W$.  Though modules and quasimodules share the same notation, the meaning of $W$ should be clear from the context. 

\begin{rem}
Any module for a M\"{o}bius vertex algebra $V$ is a quasimodule for $V$ with $f(x_1,x_2)=1$.
\end{rem}

\begin{rem}
Because, the quasi-Jacobi identity continues to hold if multiplied by a monomial of the form $cx_1^i x_2^j$ with $i,j \in \Z$ and $c\in\C$, we can assume any polynomial $f$ for $u,v\in V$ in the quasi-Jacobi identify has the form: $\sum_{L \geq i,j \geq 0} a_{ij} x_1^ix_2^j$ with $a_{00} = 1$.  
\end{rem}

By taking a suitable residue of the quasi-Jacobi identity, we obtain a formulas for a quasi-associativity identity and a quasi-commutativity identity.  These identities are analogous to the associativity and commutativity identities for modes of a vertex algebra obtained from Borcherds's identity.  The `quasi' versions of associativity and commutativity apply only to modes acting on elements of a quasimodule.  There are no algebra versions of these `quasi'-identities, since the algebra obeys the standard Jacobi identity.  The following lemma gives the quasi-associativity identity.

\begin{lem} \label{assoclem} 
Given a M\"{o}bius vertex algebra $V$ and a quasimodule $W$, for $u,v \in V$ and $m, n \in \Z$, there exists $f(x_1,x_2)=\sum_{L \geq i,j \geq 0} a_{ij} x_1^ix_2^j$ with $a_{00}=1$ such that 
\begin{eqnarray}
\sum_{L \geq i,j \geq 0} \sum_{k \geq 0} a_{ij} \binom{i}{k} (u_{m+k}v)_{n+i+j-k} 
&=& \sum_{L \geq i,j \geq 0} \sum_{k \geq 0} (-1)^k \binom{m}{k}  a_{ij} u_{m+i-k}v_{n+j+k}\nonumber \\ 
&&- \sum_{L \geq i,j \geq 0} \sum_{k \geq 0} (-1)^{k+m} \binom{m}{k}  a_{ij} v_{m+n+j-k}u_{i+k} \label{qassocid}
\end{eqnarray}
as operators on $W$.
\end{lem}

\pf Take the $\Res_{x_0}\Res_{x_1}\Res_{x_2} x_0^m x_2^n$ of both sides of the quasi-Jacobi identity.  The residue of the right hand side yields:
\begin{eqnarray}
\lefteqn{\Res_{x_0}\Res_{x_1}\Res_{x_2} x_0^m x_2^{n-1}  \d \left(\frac{x_1- x_0}{x_2}\right)f(x_1,x_2)Y_W(Y(u,x_0)v,x_2)}\\
&=&\Res_{x_0}\Res_{x_1}\Res_{x_2} x_0^m x_2^n x_1^{-1} \d \left(\frac{x_2+ x_0}{x_1}\right)f(x_1,x_2)Y_W(Y(u,x_0)v,x_2)\\
&=&\Res_{x_0}\Res_{x_2} x_0^m x_2^n f(x_2+x_0,x_2)Y_W(Y(u,x_0)v,x_2)\\
&=&\Res_{x_0}\Res_{x_2} x_0^m \sum_{L \geq i,j \geq 0}  a_{ij} (x_2+x_0)^i x_2^{n+j}Y_W(Y(u,x_0)v,x_2)\\
&=&\sum_{L \geq i,j \geq 0} \sum_{k \geq 0} a_{ij} \binom{i}{k} (u_{m+k}v)_{n+i+j-k}. 
\end{eqnarray}

The residue of the right hand side yields:
\begin{eqnarray}
\lefteqn{\Res_{x_0}\Res_{x_1}\Res_{x_2}   x_2^{n} x_0^{m-1}\d \left(\frac{x_1 - x_2}{ x_0}\right)f(x_1,x_2)Y_W(u,x_1)Y_W(v,x_2)} \nonumber \\ 
& & - x_2^{n} x_0^{m-1} \d \left(\frac{x_2- x_1 }{ -x_0}\right)f(x_1,x_2)Y_W(v,x_2)Y_W(u,x_1)\\
&=& \Res_{x_1}\Res_{x_2}  x_2^n (x_1-x_2)^m f(x_1,x_2)Y_W(u,x_1)Y_W(v,x_2) \nonumber \\
& & - x_2^{n} (-1)^m(x_2-x_1)^{m}f(x_1,x_2)Y_W(v,x_2)Y_W(u,x_1)\\
&=& \Res_{x_1}\Res_{x_2} \sum_{L \geq i,j \geq 0} a_{ij} \sum_{k \geq 0} (-1)^k \binom{m}{k} x_1^{m+i-k} x_2^{n+j+k}  Y_W(u,x_1)Y_W(v,x_2)  \nonumber \\
& & - \sum_{L \geq i,j \geq 0} a_{ij} \sum_{k \geq 0} (-1)^{k+m} \binom{m}{k}x_2^{m+n+j-k} x_1^{i+k} Y_W(v,x_2)Y_W(u,x_1)\\
&=&\sum_{L \geq i,j \geq 0} \sum_{k \geq 0} (-1)^k \binom{m}{k}  a_{ij} u_{m+i-k}v_{n+j+k} \nonumber \\
& & - \sum_{L \geq i,j \geq 0} \sum_{k \geq 0} (-1)^{k+m} \binom{m}{k}  a_{ij} v_{m+n+j-k}u_{i+k}.
\end{eqnarray}
\qed 

Two specializations, identities (\ref{replacementid}) and (\ref{straight2}),  of this quasi-associativity identity will be used in the proof of the spanning set results.  
Making the substitution $m=-2$ and isolating $(u_{-2}v)_{n}$, we obtain:

\begin{eqnarray} \label{replacementid}
(u_{-2}v)_n
&=&-\sum_{\genfrac{}{}{0pt}{}{L \ge i,j \ge 0}{i+j \ne 0}} \sum_{k \geq 0} a_{ij} \binom{i}{k} (u_{-2+k}v)_{n+i+j-k} \nonumber \\
&&+ \sum_{L \geq i,j \geq 0} \sum_{k \geq 0} (-1)^k \binom{-2}{k}  a_{ij} u_{-2+i-k}v_{n+j+k} \nonumber \\
&&- \sum_{L \geq i,j \geq 0} \sum_{k \geq 0} (-1)^{k} \binom{-2}{k}  a_{ij} v_{-2+n+j-k}u_{i+k} 
\end{eqnarray}

This identity is used to replace general expressions with expressions involving only elements of a certain generating set, and we refer to (\ref{replacementid}) as the replacement identity.  The second specialization of Lemma \ref{assoclem} will be used to impose a difference-one condition on the spanning set.  We begin by collecting modes with repeated indices in (\ref{qassocid}).

\begin{lem} \label{straightid} 
Given a M\"{o}bius vertex algebra $V$ and a quasimodule $W$, for $u,v \in V$ and $m, n \in \Z$, there exists $f(x_1,x_2)=\sum_{L \geq i,j \geq 0} a_{ij} x_1^ix_2^j$ with $a_{00}=1$ such that 
\begin{eqnarray} 
\lefteqn{\sum_{\genfrac{}{}{0pt}{}{L \ge i,j \ge 0}{i+j\in 2\Z}} a_{ij}\left( u_{n+\frac{i+j}{2}}v_{n+\frac{i+j}{2}}+v_{n+\frac{i+j}{2}}u_{n+\frac{i+j}{2}} \right)=} \nonumber \\
&& \sum_{L \geq i,j \geq 0}   \sum_{k \geq 0} a_{ij} \binom{i}{k} (u_{-1+k}v)_{2n+1+i+j-k} \nonumber \\
&-&\sum_{L \geq i,j \geq 0} \sum_{\genfrac{}{}{0pt}{}{k \geq 0}{k \ne -n-1+\frac{i-j}{2}}}   a_{ij} u_{-1+i-k}v_{2n+1+j+k} \label{uvid} \nonumber \\
&-& \sum_{L \geq i,j \geq 0} \sum_{\genfrac{}{}{0pt}{}{k \geq 0}{k \ne n+\frac{j-i}{2}}} a_{ij} v_{2n+j-k}u_{i+k} \label{vuid}
\end{eqnarray}
as operators on $W$.
\end{lem}

\pf In Lemma \ref{assoclem} let $m=-1$ and replace $n$ with $2n+1$ to obtain:

\begin{eqnarray}
\sum_{L \geq i,j \geq 0} \sum_{k \geq 0} a_{ij} \binom{i}{k} (u_{-1+k}v)_{2n+1+i+j-k} 
&=& \sum_{L \geq i,j \geq 0} \sum_{k \geq 0}   a_{ij} u_{-1+i-k}v_{2n+1+j+k}\nonumber \\ 
&&+ \sum_{L \geq i,j \geq 0} \sum_{k \geq 0}   a_{ij} v_{2n+j-k}u_{i+k}.
\end{eqnarray}

Collect repeated modes, corresponding to $k=-n-1 +\frac{i-j}{2}$ in the first sum on the right-hand side and $k=n+\frac{j-i}{2}$ in the second sum on the right-hand side to obtain the desired identity.  
\qed

Note that on the right-hand side of (\ref{vuid}),
$$\sum_{\genfrac{}{}{0pt}{}{L \ge i.j \ge 0}{ i+j\in 2\Z}} a_{ij}\left(u_{n+\frac{i+j}{2}}v_{n+\frac{i+j}{2}}+v_{n+\frac{i+j}{2}}u_{n+\frac{i+j}{2}} \right),$$
either the term $u_{n}v_{n}$ or term $v_{n}u_{n}$ appears, but not both.  The term $u_{n}v_{n}$ corresponds to $k=-n-1$ while the term $v_{n}u_{n}$ corresponds to $k=n$, but $k$ is a non-negative integer. This means that if $n<0$ then only $u_{n}v_{n}$ appears, and if $n\ge0$ then only $v_{n}u_{n}$ appears. For $n<0$ we rewrite the identity in Lemma \ref{straightid} as

\begin{eqnarray} \label{straight2}
u_nv_n=
&-&\sum_{\genfrac{}{}{0pt}{}{L \ge i,j \ge 0}{i+j \ne 0,i+j \in 2\Z }} a_{ij}\left( u_{n+\frac{i+j}{2}}v_{n+\frac{i+j}{2}}+v_{n+\frac{i+j}{2}}u_{n+\frac{i+j}{2}} \right) \nonumber \\
&+& \sum_{L \geq i,j \geq 0}   \sum_{k \geq 0} a_{ij} \binom{i}{k} (u_{-1+k}v)_{-2n+1+i+j-k} \nonumber \\
&-&\sum_{L \geq i,j \geq 0} \sum_{\genfrac{}{}{0pt}{}{k \geq 0}{k \ne -n-1+\frac{i-j}{2}}}   a_{ij} u_{-1+i-k}v_{-2n+1+j+k} \nonumber \\
&-& \sum_{L \geq i,j \geq 0} \sum_{\genfrac{}{}{0pt}{}{k \geq 0}{k \ne n+\frac{j-i}{2}}} a_{ij} v_{-2n+j-k}u_{i+k}, 
\end{eqnarray}

and for $n \ge 0$ we can isolate $v_nu_n$ in a similar manner.  This identity is used to remove modes with repeated indices in spanning set element, imposing a difference-one condition on spanning sets.  Such identities are sometimes called `straightening' identities, and this is what (\ref{straight2}) will be referred to as.

Another important identity for the construction of Poincar\'{e}-Birkhoff-Witt-like spanning sets is a commutativity identity for operators.  We obtain the following quasi-commutativity identity by taking a suitable residue of the quasi-Jacobi identity.

\begin{lem} \label{commlem} 
Given a M\"{o}bius vertex algebra $V$ and a quasimodule $W$, for $u,v \in V$ and $m, n \in \Z$, there exists $f(x_1,x_2)=\sum_{L \geq i,j \geq 0} a_{ij} x_1^ix_2^j$ with $a_{00}=1$ such that \begin{eqnarray} \label{commid1}
\sum_{L \geq i,j \geq 0} a_{ij} [u_{m+i},v_{n+j}]
=\sum_{L \geq i,j \geq 0} \sum_{k \ge 0} a_{ij} \binom{m+i}{k}(u_{k}v)_{m+n+i+j-k}.
\end{eqnarray}
as operators on $W$.
\end{lem}

\pf Take $\Res_{x_0}\Res_{x_1}\Res_{x_2} x_1^m x_2^n$ of both sides of the quasi-Jacobi identitiy.  This residue of the left-hand side yields:

\begin{eqnarray}
\lefteqn{\Res_{x_0}\Res_{x_1}\Res_{x_2}   x_1^{m} x_2^{n} x_0^{-1}\d \left(\frac{x_1 - x_2 }{ x_0}\right)f(x_1,x_2)Y_W(u,x_1)Y_W(v,x_2)}\\ \nonumber
& & - x_1^{m} x_2^{n}x_0^{-1} \d \left(\frac{x_2- x_1 }{ -x_0}\right)f(x_1,x_2)Y_W(v,x_2)Y_W(u,x_1)\\
&=& \Res_{x_1}\Res_{x_2} x_1^{m} x_2^{n} f(x_1,x_2)[Y_W(u,x_1),Y_W(v,x_2)] \\
&=& \sum_{L \geq i,j \geq 0}a_{ij} x_1^{m+i} x_2^{n+j} [Y_W(u,x_1),Y_W(v,x_2)] \\
&=&\sum_{L \geq i,j \geq 0} a_{ij} [u_{m+i},v_{n+j}].
\end{eqnarray}

This residue of the right-hand side yields:

\begin{eqnarray}
\lefteqn{\Res_{x_0}\Res_{x_1}\Res_{x_2} x_1^{m} x_2^{n-1}  \d \left(\frac{x_1- x_0}{x_2}\right)f(x_1,x_2)Y_W(Y(u,x_0)v,x_2)}\\
&=& \Res_{x_0}\Res_{x_1}\Res_{x_2} x_1^{m-1} x_2^{n}  \d \left(\frac{x_2+ x_0}{x_1}\right)f(x_1,x_2)Y_W(Y(u,x_0)v,x_2)\\
&=& \Res_{x_0}\Res_{x_2} (x_2+x_0)^{m} x_2^{n} f(x_2+x_0,x_2)Y_W(Y(u,x_0)v,x_2)\\
&=& \Res_{x_0}\Res_{x_2} \sum_{L \geq i,j \geq 0} a_{ij} (x_2+x_0)^{m+i}x_2^{n+j} Y_W(Y(u,x_0)v,x_2)\\
&=&\sum_{L \geq i,j \geq 0} \sum_{k \ge 0} a_{ij} \binom{m+i}{k}(u_{k}v)_{m+n+i+j-k}.
\end{eqnarray}

\qed 

By manipulating (\ref{commid1}), we obtain the following formula for the commutator $[u_m,v_n]$:

\begin{eqnarray} \label{commid2} 
[u_m,v_n] = - \sum_{\genfrac{}{}{0pt}{}{L \ge i,j \ge 0}{i+j \ne 0}} a_{ij} [u_{m+i},v_{n+j}] + \sum_{L \geq i,j \geq 0} \sum_{k \ge 0} a_{ij} \binom{m+i}{k}(u_{k}v)_{m+n+i+j-k}.
\end{eqnarray}

\begin{de}
The {\it degree} of an operator $u_n$ for $u$ homogeneous is $\deg(u_n)= \wt u - n -1$.
\end{de}

\begin{rem} \label{degree_id_remark}
Note that the degree of the operator on the left-hand side of identities (\ref{replacementid}), (\ref{straight2}), and (\ref{commid2}) is greater than or equal to the degree of any operator on the right-hand sides of each equation respectively.
\end{rem}

In summary, this section presents the notion of a quasimodule and develops identities for these objects.  In particular the identities (\ref{replacementid}), (\ref{straight2}), and (\ref{commid2}) are used to prove the difference-zero and difference-one quasimodule spanning sets.  These identities are derived from the quasi-associativity identity of Lemma \ref{assoclem} and the quasi-commutativity identity of Lemma \ref{commlem}.  The identities are analogous to the associativity and and commutativity identities for modes of vertex algebras.  In the spanning set proofs the replacement identity (\ref{replacementid}) and commutator identity (\ref{commid2}) are used to adjust expressions involving modes in a useful way.  The straightening identity (\ref{straight2}) is used to impose the difference-one ordering condition.

\section{Difference-zero quasimodule spanning set}

In this section, we develop a  spanning set for quasimodules for M\"{o}bius vertex algebras with ordering restrictions analogous to those of the minimal Poincar\'{e}-Birkhoff-Witt-like spanning sets for vertex operator algebras and their modules of Karel and Li \cite{MR1700507}.  These spanning sets consist of monomials of the form
\begin{eqnarray}
u^{(1)}_{n_1} \cdots u^{(r)}_{n_r}w,
\end{eqnarray}
where $\deg(u^{(1)}_{n_1})\geq \cdots \geq \deg(u^{(r)}_{n_r})\geq 0$ and $w$ is a lowest weight vector.  The ordering restriction of this spanning set is: degrees of adjacent modes must have difference-zero or more.  The difference-zero spanning set for quasimodules of M\"{o}bius vertex algebras features a slightly different ordering restriction:  the difference of indices of adjacent modes is zero or more.

For the Karel and Li module spanning set, it is possible to rearrange the monomials so that the ordering condition becomes a difference-zero condition on the indices of the modes:  $n_1 \leq \cdots \leq n_k <T$ where $T$ is some non-negative integer that depends on the elements of $V/C_1(V)$.  This $T$ is a uniform upper bound on the indices of modes that annihilate a generating vector. If $V$ is $C_1$-cofinite then such a $T$ exists, but assuming $C_1$-cofiniteness of the algebra is not necessary to prove a difference-zero spanning set.

\begin{de}
A vector $w$ of a quasimodule $W$ for a  M\"{o}bius vertex algebra $V$ is {\it uniformly annihilated} by a set of vectors $X \subseteq V$ if there exists $T \in \N$ such that $x_nw=0$ for any $ x\in X$, $n \ge T$.  The smallest such $T$ for a given $X$ is called the {\it order of uniform annihilation} of $w$ by $X$.
\end{de}

For algebras the creation axiom ensures that the vacuum vector is uniformly annihilated by elements of the algebra with $T=0$, but modules do not have a creation axiom. If $X$ is a finite set, any vector in $W$ is uniformly annihilated by $X$ because of the lower truncation property of quasimodules.  In this section $X$ is a set of a homogeneous representatives of a basis for $V/C_1(V)$.  This means that if a M\"{o}bius vertex algebra $V$ is $C_1$-cofinite, any vector in a quasimodule will be uniformly annihilated by $X$.  However we do not need to assume $C_1$-cofiniteness to prove the existence of a difference-zero quasimodule spanning set;  we only need to assume that the quasimodule is generated by a vector uniformly annihilated by this set $X$.

\begin{rem} \label{Lirem}
Proposition 3.3 in \cite{MR1676852} states that a vertex operator algebra $V$ is spanned by elements of the form:
\begin{eqnarray}
x^{(1)}_{n_1} \cdots x^{(r)}_{n_r}\1,
\end{eqnarray}
where $x^{(i)}\in X$, and $X$ a set of homogeneous representatives of a basis for $V/C_1(V)$.  This proposition holds for the  more general M\"{o}bius vertex algebras, as only the $L(-1)$ derivative is used in the proof of Proposition 3.3.
\end{rem}

\begin{lem} \label{c1reps}
For a  M\"{o}bius vertex algebra $V$, a quasimodule module $W$ generated by $w \in W$, and $X$ a set of homogeneous representatives of a basis for $V/C_1(V)$, $W$ is spanned by elements of the form
\begin{eqnarray}
x^{(1)}_{n_1} \cdots x^{(r)}_{n_r}w
\end{eqnarray}
where $r \in \N$, $x^{(1)}, \ldots, x^{(r)} \in X$, and $n_1, \ldots , n_r \in \Z$.
\end{lem}

\pf
Since $W$ is generated by $w$, a general element of $W$ has the form
\begin{eqnarray}
u^{(1)}_{n_1} \cdots u^{(r)}_{n_r}w \label{diff0lem}
\end{eqnarray}
where $r \in \N$; $u^{(1)}, \ldots, u^{(r)} \in V$; and $n_1, \ldots , n_r \in \Z$.  Each $u^{(i)}$ in  (\ref{diff0lem}) can be replaced by a sum of vectors of the form  $x^{(1)}_{n_1} \cdots x^{(r)}_{n_r}\1$ by Remark \ref{Lirem}.  Apply the quasi-associativity identity, Lemma \ref{assoclem}.

\begin{eqnarray}
(x^{(1)}_{n_1} \cdots x^{(r)}_{n_r}\1)_{m}  
&=&-\sum_{\genfrac{}{}{0pt}{}{L \ge i,j \ge 0}{i+j \ne 0}} \sum_{k \geq 0} a_{ij} \binom{i}{k} (x^{(1)}_{n_1+k}x^{(2)}_{n_2} \cdots x^{(r)}_{n_r}\1)_{m+i+j-k}  \label{quickref1}\\
&&+ \sum_{L \geq i,j \geq 0} \sum_{k \geq 0} (-1)^k \binom{n_1}{k}  a_{ij} x^{(1)}_{-2+i-k}(x^{(2)}_{n_2} \cdots x^{(r)}_{n_r}\1)_{m+j+k}  \\
&&- \sum_{L \geq i,j \geq 0} \sum_{k \geq 0} (-1)^{k+n_1} \binom{n_1}{k}  a_{ij} (x^{(2)}_{n_2} \cdots x^{(r)}_{n_r}\1)_{n_1+m+j-k}x^{(1)}_{i+k} 
\end{eqnarray}

Since $i+j>0$ in (\ref{quickref1}), 
\begin{eqnarray}
\deg (x^{(1)}_{n_1+k}x^{(2)}_{n_2} \cdots x^{(r)}_{n_r}\1)_{m+i+j-k}>\deg (x^{(1)}_{n_1} \cdots x^{(r)}_{n_r}\1)_m
\end{eqnarray}
for all $k$.  This means than the process of repeatedly applying Lemma \ref{assoclem} will eventually terminate because of the lower truncation property of $W$.  All we are left with is a sum of monomials of terms of the form $x_n$ with $x\in X$ and $n \in \Z$.
\qed

Given an $\N$-graded M\"{o}bius vertex algebra $V$ and a quasimodule $W$, we define a filtration on $W$:
\begin{eqnarray}
W^{(0)} \subset  W^{(1)} \subset \cdots \subset  W^{(s)} \subset \cdots \subset  W,
\end{eqnarray}
where $W^{(s)}=\span\{ u^{(1)}_{n_1} \cdots u^{(r)}_{n_r}w : \sum_{i=1}^r \wt u^{(i)} \le s\}$ for homogenous $u^{(i)} \in V$.  A crucial feature of this filtration, proved below, is that rearrangement of modes does not change the filtration level of a monomial.  More specifically two monomials consisting of the same modes, but with different orderings, are equal up to addition by an element of a lower filtration level.

\begin{lem}\label{fillem1}
Let $W$ be a quasimodule for a M\"{o}bius vertex algebra. For a monomial $u^{(1)}_{n_1} \cdots u^{(r)}_{n_r}w \in W^{(s)}$,

\begin{eqnarray}u^{(1)}_{n_1} \cdots u^{(r)}_{n_r}w=u^{(\sigma(1))}_{n_{\sigma(1)}} \cdots u^{(\sigma(r))}_{n_{\sigma(r)}}w + R,\end{eqnarray}

where $\sigma \in Sym(r)$ and $R \in W^{(s-1)}$.
\end{lem}

\pf It is sufficient to prove that this is true for a transposition of adjacent modes.  Furthermore, it is sufficient to prove for any $u,v \in V$, $w \in W$, and $n,m \in \Z$ that $u_nv_mw=v_mu_nw + R$ where $u_nv_mw  \in W^{(s)}$ and $R \in W^{(s-1)}$.   Since $u_nv_mw=v_mu_nw + [u_n,v_m]w$, it is sufficient to show that $[u_m,v_n]w \in W^{(s-1)}$.  We use the commutator identity (\ref{commid2}):

\begin{eqnarray}
[u_m,v_n]w = - \sum_{\genfrac{}{}{0pt}{}{L \ge i,j \ge 0}{i+j \ne 0}} a_{ij} [u_{m+i},v_{n+j}]w + \sum_{L \geq i,j \geq 0} \sum_{k \ge 0} a_{ij} \binom{m+i}{k}(u_{k}v)_{m+n+i+j-k}.
\end{eqnarray}

Because $\wt (u_{k}v)= \wt u +\wt v -k-1<\wt u +\wt v$, all the terms in the second (triple) sum are in $W^{(s-1)}$.  We now argue that elements in the first (single) sum can be rewritten as elements of a lower filtration level by a finite process. In the terms of the form $[u_{m+i},v_{n+j}]$, $i+j>0$, so the $\deg([u_{m+i},v_{n+j}]) > \deg([u_m,v_n])$.  We again apply (\ref{commid2}) to all these commutator terms  of the form $[u_{m+i},v_{n+j}]$ in the first (single) sum.  The result of applying (\ref{commid2}) is a sum commutator elements with strictly larger degrees than in the previous step and other terms of a strictly lower filtration level.   The lower truncation property ensures that for large enough L, $u_Lw=0$ and $v_Lw=0$.  We can continue to apply identity (\ref{commid2}) until all we are left with are terms of the form $(u_{p}v)_{q}$ with $p \ge 0$ which are in $W^{(s-1)}$. \qed

We now prove a difference-zero spanning set element for quasimodules for $\N$-graded M\"{o}bius vertex algebra.  The assumption of uniform annihilation by $X$ is always satisfied if $V$ is $C_1$-cofinite where $X$ is a set of representatives of a basis for $V/C_1(V)$.

\begin{thm}\label{diff0thm}
For an $\N$-graded M\"{o}bius vertex algebra $V$, $X$ a set of homogeneous representatives of a basis for $V/C_1(V)$, and a quasimodule module $W$ generated by a vector $w$ that is uniformly annihilated by $X$, $W$ is spanned by the elements of the form
\begin{eqnarray}
x^{(1)}_{n_1} \cdots x^{(r)}_{n_r}w
\end{eqnarray}
with $n_1\leq \cdots \leq n_r<T$ where $r \in \N$, $x^{(1)}, \ldots, x^{(r)} \in X$, $n_1, \ldots , n_r \in \Z$, and $T$ is the order of uniform annihilation of $w$ by $X$.
\end{thm}

\pf
We prove this theorem by induction on the filtration level of the module $W$.  By Lemma \ref{c1reps} an element in $W^{(0)}$ can be written a linear combination of vectors the form $x^{(1)}_{n_1} \cdots x^{(r)}_{n_r}w$ with the $x_i$'s in $X$ with weight $0$.  Applying Lemma \ref{fillem1} to these vectors, we may rearrange them so that the modes have the desired ordering. Since $x^{(1)}_{n_1} \cdots x^{(r)}_{n_r}w \in W^{(0)}$, after rearrangement the remainder term $R$ is  $0$.    

Now assuming the induction hypothesis holds for lower filtration levels and given an element of $W^{(s)}$, we may write this element in terms of elements of $X$ by Lemma \ref{c1reps}.  By Lemma \ref{fillem1} we may rearrange the modes to display the desired order and apply the induction hypothesis to the remainder term.\qed

Using the commutator identity  we can impose a different ordering on the modes of quasimodules. The following is the ordering of Karel and Li in their difference-zero spanning set for vertex operator algebras and modules.

\begin{cor}
For an $\N$-graded M\"{o}bius vertex algebra $V$, $X$ a set of homogeneous representatives of a basis for $V/C_1(V)$, and a quasimodule module $W$ generated by a vector $w$ that is uniformly annihilated by $X$, $W$ is spanned by the elements of the form
\begin{eqnarray}
x^{(1)}_{n_1} \cdots x^{(r)}_{n_r}w
\end{eqnarray}
with $\deg (x^{(1)}_{n_1}) \geq \cdots \geq \deg (x^{(r)}_{n_r}) \geq -T-1$ where $r \in \N$, $x^{(1)}, \ldots, x^{(r)} \in X$, $n_1, \ldots , n_r \in \Z$, and $T$ is the order of uniform annihilation of $w$ by $X$.
\end{cor}

\pf This follows from the same proof of \ref{diff0thm}.  Use (\ref{commid2}) to impose the desired ordering. \qed




This difference-zero spanning set for quasimodules of M\"{o}bius vertex algebras is generated by representatives of a basis for $V/C_1(V)$.  In next section, the difference-one spanning set for quasimodules is generated by a larger set, representatives of a basis for $V/C_2(V)$.  For the difference-zero spanning set, the tradeoff is a smaller generating set for less restrictive ordering restrictions ordering restrictions.  For the difference-one spanning set, the trade off is the opposite.

\section{Difference-one quasimodule spanning set}

In this section, we prove a difference-one spanning set for quasimodules for $\N$-graded M\"{o}bius vertex algebras.  This is a generalization to quasimodules of similar difference-one spanning sets for vertex operator algebras of CFT type \cite{MR1990879}, modules for vertex operator algebras \cite{MR1927435}, vertex operator algebras with a non-negative grading \cite{MR2046807}, and twisted modules for vertex operator algebras \cite{MR2039213}.

Given an $\N$-graded M\"{o}bius vertex algebra $V$ and a quasimodule $W$ generated by $w \in W$,  the quasimodule $W$ is spanned by elements of the form $u^{(1)}_{n_1} \cdots u^{(r)}_{n_r}w$ where $u^{(i)} \in V$ and $n_i \in \Z$. We demonstrate a spanning set that resricts such expressions in two ways:  the vectors $u^{(i)}$ will be limited to a finite subset of $V$, and the indices $n_i $ will be strictly increasing.

We again use the filtration for a $V$-quasimodule $W$ introduced in the previous section:
\begin{eqnarray}
W^{(0)} \subset  W^{(1)} \subset \cdots \subset  W^{(s)} \subset \cdots \subset  W,
\end{eqnarray}
where $W^{(s)}=\span\{ u^{(1)}_{n_1} \cdots u^{(r)}_{n_r}w : \sum_{i=1}^r \wt u^{(i)} \le s\}$ for homogenous $u^{(i)} \in V$. This filtration has two important properties: any rearrangement of the the modes in a monomial $u^{(1)}_{n_1} \cdots u^{(r)}_{n_r}w$ does not alter the filtration level, and  any replacement of a vector $u^{(i)}$  by its representative in $V/C_2(V)$ does not alter the filtration level.  Lemma \ref{fillem1} verifies this first property and the following lemma verifies this second property.


\begin{lem} \label{fillemc2}
Let $W$ be a quasimodule for a M\"{o}bius vertex algebra. For a monomial $u^{(1)}_{n_1} \cdots u^{(r)}_{n_r}w \in W^{(s)}$,

\begin{eqnarray}u^{(1)}_{n_1} \cdots u^{(r)}_{n_r}w=x^{(1)}_{n_{1}} \cdots x^{(r)}_{n_{r}}w + R\end{eqnarray}

where $x^{(1)}_{n_{1}} \cdots x^{(r)}_{n_{r}}w \in W^{(s)}$, $R \in W^{(s-1)}$, and  $x^{(i)}$ is a representative of $u^{(i)} + C_2(V)$ for $1 \ge i \ge r$ .
\end{lem}

\pf  By linearity, it is sufficient to show  that $y_nw=x_nw+(u_{-2}v)_nw$  with $y_nw, x_nw \in W^{(s)}$ and $(u_{-2}v)_nw \in W^{(s-1)}$ for any homogeneous $x$ and $y$ with $y=x +u_{-2}v$.  We use the replacement identity (\ref{replacementid}).

\begin{eqnarray}
(u_{-2}v)_nw
&=&-\sum_{\genfrac{}{}{0pt}{}{L \ge i,j \ge 0}{i+j \ne 0}} \sum_{k \geq 0} a_{ij} \binom{i}{k} (u_{-2+k}v)_{n+i+j-k}w \label{ass1}\\
&&+ \sum_{L \geq i,j \geq 0} \sum_{k \geq 0} (-1)^k \binom{-2}{k}  a_{ij} u_{-2+i-k}v_{n+j+k}w \label{ass2}\\
&&- \sum_{L \geq i,j \geq 0} \sum_{k \geq 0} (-1)^k \binom{-2}{k}  a_{ij} v_{-2+n+j-k}u_{i+k}w \label{ass3}
\end{eqnarray}

If $y_nw \in W^{(s)}$, then $s=\wt y=\wt u_{-2}v=\wt u +\wt v +1$.  So $\wt u +\wt v=s-1$  and all the terms in (\ref {ass2})and (\ref{ass3}) are in $W^{(s-1)}$.  Since $\wt(u_{-2+k}v)= \wt u +\wt v -k+1=s-k$, all the  terms in (\ref{ass1}) are in $W^{(s-1)}$ except for 
\begin{eqnarray}
\sum_{\genfrac{}{}{0pt}{}{L \ge i,j \ge 0}{i+j \ne 0}} a_{ij} (u_{-2}v)_{n+i+j}w.
\end{eqnarray}
In this sum at least one of $i$ or $j$ is positive, and $n+i+j>n$.  We  apply the replacement identity (\ref{replacementid}) to the elements of  the form $(u_{-2}v)_{n+i+j}w$.  This results in  elements in $W^{(s-1)}$ and a sum elements of the form $(u_{-2}v)_{K}w$ with $K>n+i+j$ . Repeating this process of applying the replacement identity to elements of the form $(u_{-2}v)_{K}w$ allows $K$ to grow arbitrarily large.  For large enough $K$, $(u_{-2}v)_{K}w=0$, leaving only elements in $W^{(s-1)}$. \qed


We now prove the main result, a difference-one spanning set for quasimodules.  This proof is an induction argument on module elements of a certain filtration level and certain length.  Throughout the proof, we make use of the replacement and reordering properties of the module filtration.  This proof is similar in spirit to Yamauchi's proof of a order spanning set for twisted modules \cite{MR2039213}.

\begin{thm}\label{mainthm}
For an $\N$-graded M\"{o}bius vertex algebra $V$, $X$ a set of homogeneous representatives of a basis for $V/C_2(V)$, and a quasimodule module $W$ generated by a vector $w$ that is uniformly annihilated by $X$, $W$ is spanned by the elements of the form
\begin{eqnarray}
x^{(1)}_{n_1} \cdots x^{(r)}_{n_r}w
\end{eqnarray}
with $n_1< \cdots < n_r<T$ where $r \in \N$; $x^{(1)}, \ldots, x^{(r)} \in X$; $n_1, \ldots , n_r \in \Z$; and $T$ is the order of uniform annihilation of $w$ by $X$.
\end{thm}

\pf   The quasimodule $W$ is generated by $w$, and $W$ is spanned by elements of the form:
\begin{eqnarray}
u^{(1)}_{n_1} \cdots u^{(r)}_{n_r}w
\end{eqnarray}
with $r \in \N$; $u^{(1)}, \ldots, u^{(r)} \in V$; and $n_1, \ldots , n_r \in \Z$.  We procede by induction on pairs $(s,r) \in \N \times \N$ with the ordering $(s_1,r_1)>(s_2,r_2)$ if $s_1>s_2$ or $s_1=s_2$ and $r_1>r_2$, where $s$ is the filtration level of a monomial $u^{(1)}_{n_1} \cdots u^{(r)}_{n_r}w$, and $r$ is its length. 

For $(s,0)$ with $s\in \N$, the only element is $w$, which is in the desired form.  
Now we assume the induction hypothesis: for all $(s',r')<(s,r)$, any vector of induction level $(s',r')$ can be rewritten as a linear combination of vectors of the desired form: $x^{(1)}_{n_1} \cdots x^{(l)}_{n_l}w$
where $l \in \N$; $x^{(1)}, \ldots, x^{(l)} \in X$; $n_1, \ldots , n_l \in \Z$; $n_1< \cdots < n_l<T$; $l \le r$; and $x^{(1)}_{n_1} \cdots x^{(l)}_{n_l}w \in W^{(s')}$.  That is, the induction hypothesis is: any vector of filtration level $s'$ or and vector of filtration level $s$ and of shorter length can be rewritten as a linear combination of monomials of the desired form.

Now consider a monomial $u^{(1)}_{n_1} \cdots u^{(r)}_{n_{r}}w \in W^{(s)}$.  By applying the `rearrangement' and `replacement lemmas', Lemma \ref{fillem1} and Lemma \ref{fillemc2}, we rewrite this monomial as
\begin{eqnarray}x^{(1)}_{n_1} \cdots x^{(r)}_{n_{r}}w +R\end{eqnarray}
with $x^{(1)}, \ldots, x^{(r)} \in X$; $n_1 \leq \cdots \leq n_r<T$; and $R \in W^{(s-1)}$.  By the induction hypothesis $R$ can be rewritten without repetition of modes.  We must show that $x^{(1)}_{n_1} \cdots x^{(r)}_{n_{r}}w$ can be rewritten without repetition of modes.  

We apply the inductive hypothesis to $x^{(2)}_{n_2} \cdots x^{(r)}_{n_{r}}w$ to obtain
\begin{eqnarray}
x^{(1)}_{n_1} \cdots x^{(r)}_{n_{r}}w=\sum x^{(1)}_{n_1} y^{(2)}_{m_2} \cdots y^{(r)}_{m_{r}}w+R
\end{eqnarray}
with $m_2< \cdots < m_r<T$ and $R$ a sum of vectors in a lower filtration.  The index of summation is suppressed for clarity.  The induction hypothesis is applied to $R$.  This leaves vectors  of the form $x^{(1)}_{n_1} x^{(2)}_{n_2} \cdots x^{(r)}_{n_{r}}w$ with $n_2< \cdots < n_r<T$  and $n_1$'s relation to $n_2$ unknown.

If  $n_1 < m_2$ the vector has the desired form, so we assume $n_1 \ge m_2$.  We must show that vectors of the form $x^{(1)}_{n_1} \cdots x^{(r)}_{n_{r}}w$ with $n_1 \geq n_2< \cdots < n_r<T$ can be rewritten a sum of vectors of the desired form.  In particular, this means me must rewrite these vectors to satisfy a difference-one condition on the indices of modes.

If  $n_1 > n_2$,  consider the following process. We use the commutator identity (\ref{commid2}) to transpose the first two modes.
\begin{eqnarray} \label{diff2less}
x^{(1)}_{n_1}x^{(2)}_{n_2} \cdots x^{(r)}_{n_{r}}w=x^{(2)}_{n_2}x^{(1)}_{n_1} x^{(3)}_{n_3} \cdots x^{(r)}_{n_{r}}w +R
\end{eqnarray}
Again we can apply the induction hypothesis to $R$. Now apply the induction hypothesis to $x^{(1)}_{n_1} x^{(3)}_{n_3} \cdots x^{(r)}_{n_{r}}w$ to get the vectors of the form
\begin{eqnarray}
y^{(2)}_{m_2} \cdots y^{(r)}_{m_{r}}w
\end{eqnarray}
with $n_2< \cdots < n_r<T$. Replace $x^{(2)}_{n_2}x^{(1)}_{n_1} x^{(3)}_{n_3} \cdots x^{(r)}_{n_{r}}w$ with $y^{(2)}_{m_2} \cdots y^{(r)}_{m_{r}}w$.  This produces the vector
\begin{eqnarray} \label{diff1step}
y^{(1)}_{m_1} \cdots y^{(r)}_{m_{r}}w
\end{eqnarray}
with $y^{(1)}= x^{(2)}$, ${m_1}={n_2}$, $n_2< \cdots < n_r<T$,and most importantly $m_1< n_1$.  This process of rearrangement and application of the induction hypothesis results in a monomial of the same filtration level, but with a strictly smaller first index.  We compare the indices of the first two modes ${m_1}$ and ${m_2}$ in \ref{diff1step}.  If  $m_1 < m_2$, we have a vector in the desired form. If $m_1 > m_2$, we repeat this process again.  Each time we repeat this process the index of leading mode $m_1$ decreases, the length $r$ of the vector is unchanged, and the degree of the operator $y^{(1)}_{m_1} \cdots y^{(r)}_{m_{r}}$ remains the same.  

Let $k=\deg \left(y^{(1)}_{m_1} \cdots y^{(r)}_{m_{r}} \right)$. So
\begin{eqnarray}
\sum^r_{i=1} (\wt y_i -m_i-1)=k,
\end{eqnarray}
and 
\begin{eqnarray}
\sum^r_{i=1}m_i&=&\sum^r_{i=1} (\wt y_i -1)-k\\
& \geq& -(r+k)
\end{eqnarray}
since $\wt y^{(i)} \geq 0$.  As we repeat this process, $m_1$ decreases, the other $m_i$'s must increase because the left had sum is bounded below by a fixed number.  For small enough $m_1$, the index of the last mode $m_r$ will be large enough so that  $y^{(r)}_{m_{r}}w=0$.  This ensures that this iterative process terminates after a finite number of steps, leaving a vector with the desired ordering.  If at any point $m_1=m_2$, we use the following process that results in monomials with $m_1 \ne m_2$.

Finally, we are left with the case where $n_1=n_2< \cdots <n_3$.  We use the straightening identity (\ref{straight2}).  We assume that $n_1<0$ which means there is no $x^{(2)}_{n_1}x^{(1)}_{n_1}$ term on the left hand side in identity (\ref{straight2}).  The argument for $n_1 \ge 0$ is similar.  

\begin{eqnarray}
\lefteqn{x^{(1)}_{n_1}x^{(2)}_{n_1} x^{(3)}_{n_3}\cdots x^{(r)}_{n_{r}}w=} \nonumber \\
&-&\sum_{\genfrac{}{}{0pt}{}{L \ge i.j > 0}{i+j \ne 0, i+j\in 2\Z}} a_{ij}\left( x^{(1)}_{n_1+\frac{i+j}{2}}x^{(2)}_{n_1+\frac{i+j}{2}}x^{(3)}_{n_3}\cdots x^{(r)}_{n_{r}}w+x^{(2)}_{n+\frac{i+j}{2}}x^{(1)}_{n+\frac{i+j}{2}} x^{(3)}_{n_3}\cdots x^{(r)}_{n_{r}}w\right) \label{line1}\\
&+& \sum_{L \geq i,j \geq 0}   \sum_{k \geq 0} a_{ij} \binom{i}{k} \left(x^{(1)}_{-1+k}x^{(2)}\right)_{2n_1+1+i+j-k}x^{(3)}_{n_3}\cdots x^{(r)}_{n_{r}}w \label{line2} \\
&-&\sum_{L \geq i,j \geq 0} \sum_{\genfrac{}{}{0pt}{}{k \geq 0}{k \ne -n-1+\frac{i-j}{2}}}   a_{ij} x^{(1)}_{-1+i-k}x^{(2)}_{2n_1+1+j+k}x^{(3)}_{n_3}\cdots x^{(r)}_{n_{r}}w \label{line3}\\
&+& \sum_{L \geq i,j \geq 0} \sum_{\genfrac{}{}{0pt}{}{k \geq 0}{k \ne n_1+\frac{j-i}{2}}} a_{ij} x^{(2)}_{2n_1+j-k}x^{(1)}_{i+k} x^{(3)}_{n_3}\cdots x^{(r)}_{n_{r}}w \label{line4}
\end{eqnarray}

Each monomial in (\ref{line2}) is in a lower filtration level because each monomial has length $r-1$, so we can apply the induction hypothesis.  In (\ref{line3}) either $-1+i-k$ or $2n_1+1+j+k$ is strictly less than $n_1$, and in (\ref{line4}) $2n_1+j-k$ or $i+k$ is strictly less than $n_1$.  So for monomials in (\ref{line3}) and (\ref{line4}), after rearrangement of modes we return to the case of monomials of the form:
\begin{eqnarray}v^{(1)}_{m_1} v^{(2)}_{m_2} \cdots v^{(r)}_{m_{r}}w\end{eqnarray}
with $m_1<n_1$ and $m_2< \cdots <m_r$.

We are left to account for the monomials in (\ref{line1}).  All of these have repeated modes with $n_1+\frac{i+j}{2}>n_1$.  We can reapply Lemma (\ref{straightid}) with $n=n_1+\frac{i+j}{2}$ and proceed as above.  Again this process eventually terminates, since we are increasing the indices with each application of the straightening identity (\ref{straight2}).  This verifies the induction hypothesis for vectors in $W^{(s)}$ proving the theorem. \qed


One of the implications of the Poincar\'{e}-Birkhoff-Witt-like spanning sets for modules is: `size' restrictions on the algebras imply that modules are similarly restricted.  The following result extends this implication to quasimodules.

\begin{cor}\label{cor1}
If $W$ is a quasimodule for a $C_2$-cofinite $\N$-graded M\"{o}bius vertex algebra and $W$ is generated by a vector $w \in W$, then $W$ is $C_n$-cofinite for $n\ge 2$ and is $c_1$-finite.
\end{cor}

\pf Since $V$ is $C_2$-cofinite $w$ is uniformly annihilated by a set $X$ of representatives of a basis of $V/C_2(V)$.  By Theorem \ref{mainthm} $W$ is spanned by elements of the form
\begin{eqnarray}
x^{(1)}_{n_1} \cdots x^{(r)}_{n_r}w
\end{eqnarray}
with $n_1< \cdots < n_r<T$ where $r \in \N$; $x^{(1)}, \ldots, x^{(r)} \in X$; $n_1, \ldots , n_r \in \Z$; and $T$ is order of uniform annihilation of $w$ by $X$.  Since the indices of the modes are strictly increasing, the maximum length of a spanning set element not in $C_n(W)$ is $T+n-1$.  Since $X$ is finite and the length of spanning set elements not in $C_n(W)$ is finite, the spanning set elements not in $C_n(W)$ is finite, and $W$ is $C_n$-cofinite.  A similar argument holds for $c_1$-cofiniteness.\qed

This argument holds if $W$ is finitely generated as well.  So any finitely generated quasimodule for a $C_2$-cofinite $\N$-graded M\"{o}bius vertex algebra is $C_n$-cofinite for $n\ge 2$ and is $c_1$-finite.  For M\"{o}bius vertex algebras, $C_2$-cofiniteness is equivalent to $C_n$-cofiniteness for $n \geq 2$.  The following corollary extends this result to quasimodules for $C_2$-cofinite $\N$-graded M\"{o}bius vertex algebras.

\begin{cor}
For a quasimodule $W$ for $C_2$-cofinite $\N$-graded M\"{o}bius vertex algebra, $C_2$-cofiniteness of $W$ is equivalent to $C_n$-cofiniteness of $W$ for $n \geq 2$.
\end{cor}

\pf The $L(-1)$-derivative property for quasimodules ensures that $C_2(W) \subseteq C_n(W)$, and $C_n$-cofiniteness of $W$ implies $C_2$-cofiniteness of $W$.  Corollary \ref{cor1} implies the converse.  \qed

  



\end{document}